\documentclass{elsarticle}
\usepackage{amsmath}
\usepackage{amssymb}
\usepackage{amsfonts}
\usepackage{amsthm}
\usepackage{rotating}
\usepackage{graphicx}
\usepackage{floatflt,epsfig}
\usepackage{lineno,hyperref}
\usepackage{enumerate}
\usepackage{colortbl}
\usepackage{array,tabularx,tabulary,booktabs}
\usepackage{longtable}
\usepackage{multirow}
\usepackage{wrapfig}
\usepackage{subcaption}
\newcolumntype{^}{>{\currentrowstyle}}

\journal{arXiv}

\def\ms{\medskip}

\newtheorem{theorem}{Theorem}
\newtheorem{lemma}{Lemma}
\newtheorem{corollary}{Corollary}

\newcommand{\VEC}[3]{{#1}_{#2},\ldots,{#1}_{#3}}
\newcommand{\floor}[1]{\left\lfloor#1\right\rfloor}
\newcommand{\ceil}[1]{\left\lceil#1\right\rceil}
\newcommand{\FFR}[2]{\floor{\frac{#1}{#2}}}
\newcommand{\CFR}[2]{\ceil{\frac{#1}{#2}}}
\newcommand{\SE}[3]{\sum\limits_{{#1}={#2}}^{#3}}

\let \CH=\binom

\begin{document}
\renewcommand{\abstractname}{Abstract}
\renewcommand{\refname}{References}
\renewcommand{\tablename}{Table.}
\renewcommand{\arraystretch}{0.9}
\thispagestyle{empty}
\sloppy

\begin{frontmatter}
\title{Nonexistence of certain edge-girth-regular graphs}

\author{Leen Droogendijk}
\ead{drooge001@kpnmail.nl}

\begin{abstract}
	Edge-girth-regular graphs (abbreviated as \emph{egr} graphs) are regular graphs in which every edge is contained in the same number of shortest cycles. We prove that there is no $3$-regular \emph{egr} graph with girth $7$ such that every edge is on exactly $6$ shortest cycles, and there is no $3$-regular \emph{egr} graph with girth $8$ such that every edge is on exactly $14$ shortest cycles. This was conjectured by Goedgebeur and Jooken~\cite{GJ1}.
	A few other unresolved cases are settled as well.
\end{abstract}

\begin{keyword} Extremal problems; Edge-girth-regular graphs
\vspace{\baselineskip}
\end{keyword}
\end{frontmatter}

\section{Introduction}

In 2018, Jajcay, Kiss and Miklavi\v{c}~\cite{JKM18} introduced a new type of regularity called edge-girth-regularity. They define for integers $v$, $k$, $g$ and $\lambda$ an edge-girth-regular $(v,k,g,\lambda)$ graph (abbreviated as an $egr(v,k,g,\lambda)$ graph) as a $k$-regular graph with girth $g$ on $v$ vertices such that every edge is contained in exactly $\lambda$ cycles of length $g$.
Goedgebeur and Jooken~\cite{GJ1} did an exhaustive search to find certain small \emph{egr} graphs and conjectured that there are no $egr(v,3,7,6)$-, $egr(v,3,8,10)$-, $egr(v,3,8,12)$-, and $egr(v,3,8,14)$-graphs.
This paper proves that conjecture to be true.
The existence of $egr(v,4,4,4)$-graphs was unresolved as well (see \cite{DFJR1}), and we show that they do not exist either.
Finally the unresolved cases $egr(v,4,5,8)$ and $egr(v,6,5,24)$ are shown not to exist.

\section{Some notational conventions.}\label{definitions}

We will write $P_t$ for the path graph on $t$ vertices, and $C_t$ for the cycle graph on $t$ vertices.

For a vertex $x$ of a graph $G$ we define $D_d(x)$ (or simply $D_d$ if $x$ is clear) to be the induced subgraph of all vertices at distance at most $d$ from $x$, and $B_d(x)$ (or $B_d$) to be the induced subgraph of all vertices at distance exactly $d$ from $x$.

Throughout the paper, $G$ will be some $egr(v,k,g,\lambda)$-graph, where
\begin{enumerate}
	\item $v$ is the number of vertices, a positive integer.
	\item $k$ is the degree, a positive integer, at least $3$.
	\item $g$ is the girth, a positive integer, at least $3$.
	\item $\lambda$ is the number of shortest cycles per edge.
\end{enumerate}

\section{A few lemmas.}

This section states a few lemmas that will be used in the other sections.
Not all cases need all the lemmas.
If you are interested in a special case, just skip this section, and only return to it if your case needs parts of this.

\begin{lemma}\label{egr: counting cycles using a vertex}
	Let $G$ be an $egr(v,k,g,\lambda)$-graph.
	Then every vertex of $G$ is contained in exactly $\frac{k\lambda}2$ shortest cycles.
	\begin{proof}
		Each of the $k$ edges incident with a given vertex $u$ is contained in $\lambda$ shortest cycles, and we count every cycle twice.
	\end{proof}
\end{lemma}

\begin{lemma}\label{egr: forbid certain egr subgraphs}
	Let $G$ be an $egr(v,k,g,\lambda)$-graph.
	$G$ does not contain any $egr(w,k-1,g,\lambda)$-subgraph.
	\begin{proof}
		Suppose $H$ is such subgraph, and $u$ a vertex of $H$.
		Since the degree of $H$ is only $k-1$, there is a unique edge $e$ with endpoint $u$, that is not an edge of $H$.
		Let $C$ be any shortest cycle containing $e$, then $e$ must also contain an edge $uw$ of $H$, since $u$ is on $C$ and $e$ is the only edge incident with $u$ that is not in $H$.
		But this brings the number of shortest cycles containing $uw$ to $\lambda+1$. Contradiction.
	\end{proof}
\end{lemma}

\begin{lemma}\label{egr intersections of shortest cycles}
	Let $G$ be an $egr(v,k,g,\lambda)$-graph.
	If $2$ shortest cycles of $G$ intersect, and the intersection contains at least one edge, then the intersection is connected.
	If $g$ is odd, the intersection of $2$ shortest cycles of $G$ is connected (or empty).
	\begin{proof}
		Suppose the claim is not true, so we can find shortest cycles $H_1$ and $H_2$ that intersect in a subgraph $I$ that contains an edge $pq$ and a vertex $x$ that is in a different connected component of $I$ on $H_1$.
		We may assume there is a $p,x$-path $P$ on $H_1$ of length less than $\frac g2$ (if not we exchange $p$ and $q$).
		There also is a $p,x$-path $Q$ on $C_2$ of length at most $\frac g2$.
		The paths $P$ and $Q$ combine to a closed walk of length less than $g$, so they contain a cycle of length less than $g$. Contradiction.
		
		For odd $g$ there is always a path of length less than $\frac g2$ between $2$ vertices on a shortest cycle, so we can slightly modify above argument to prove the second claim.
	\end{proof}
\end{lemma}

\begin{lemma}\label{egr: intersection of cycles}
	Let $G$ be an $egr(v,k,g,\lambda)$-graph, and let $u=\CFR{g+3}2$.
	Then a subgraph $P_u$ of $G$ is contained in at most one shortest cycle.
	\begin{proof}
		Suppose $P$ is a $P_u$ that is contained in $2$ shortest cycles. Let $v$ and $w$ be the endpoints of $P$.
		The path graph $P_u$ has length $\CFR{g+1}2$, so the $2$ $v,w$-paths completing the cycles combine to a closed walk of length $2(g-\CFR{g+1}2)=2\FFR{g-1}2<g$, which contains a cycle of length less than $g$. Contradiction.
	\end{proof}
\end{lemma}

\begin{lemma}\label{egr: intersection of cycles 2}
	Let $G$ be an $egr(v,k,g,\lambda)$-graph, and let $u=\CFR{g+3}2$.
	Then a subgraph $P_{u-i}$ in $G$ is contained in at most $(k-1)^i$ shortest cycles for $i\in\{0,\ldots,u-3\}$.
	\begin{proof}
		We use induction on $i$. The induction base, $i=0$, is Lemma \ref{egr: intersection of cycles}.
		Suppose $i>0$ and subgraph $P=P_{u-i}$ is on $(k-1)^i+1$ shortest cycles. Let $p$ be an endpoint of $P$, and let $\VEC v1{k-1}$ be the neighbors of $p$ that are not on $P$.
		By the pigeonhole principle some edge $pv_i$ must be contained in at least $(k-1)^{i-1}+1$ of the cycles, violating the induction hypothesis.
	\end{proof}
\end{lemma}

\begin{lemma}\label{egr: non-incident edges}
	Let $G$ be an $egr(v,k,g,\lambda)$-graph, with $g\geq5$.
	Then $2$ non-incident edges of $G$ can be on at most $(k-1)^{\CFR{g-5}2}$ shortest cycles.
	\begin{proof}
		Suppose $e$ and $f$ are non-incident edges of $G$, and suppose they are on $2$ shortest cycles.
		By Lemma \ref{egr intersections of shortest cycles} the intersection must be a path $P$ containing both $e$ and $f$.
		Let $Q$ be the shortest subpath of $P$ connecting $e$ and $f$, so $Q$ has length at least $1$.
		
		By Lemma \ref{egr: intersection of cycles} $P$ is at most a $P_{u-1}$ with $u=\CFR{g+3}2$, so $P$ has length at most $u-2=\CFR{g-1}2$. This means that $Q$ has length at most $\CFR{g-5}2$.
		
		By Lemma \ref{egr: intersection of cycles 2} (taking $u-i=4$, so $i=u-4=\CFR{g-5}2$), any $P_4$ is on at most $(k-1)^{\CFR{g-5}2}$ shortest cycles, so certainly $P$ is on at most $(k-1)^{\CFR{g-5}2}$ shortest cycles.
		
		If any other cycle contains both $e$ and $f$, we find a second path $Q'$ of length at most $\CFR{g-5}2$ connecting $e$ and $f$.
		Then $Q$ and $Q'$, together with $e$ and $f$ contain a cycle of length at most $2\CFR{g-5}2+2\leq2(\frac{g-5}2+1)+2=g-1$. Contradiction.
	\end{proof}
\end{lemma}

\begin{corollary}\label{egr: connectivity 2}
	Let $G$ be a connected $egr(v,k,g,\lambda)$-graph, with $g\geq5$, $t$ a positive integer.
	When $\lambda>(t-1)(k-1)^{\CFR{g-5}2}$, removal of $t$ pairwise non-incident edges cannot disconnect $G$.
	(This almost says that $G$ is ($t+1$)-edge-connected, but of course removal of $k$ edges that are all incident to the same vertex \emph{does} disconnect $G$.)
	\begin{proof}
		Suppose $\VEC e1t$ are pairwise non-incident edges that disconnect $G$.
		Each of the $\lambda$ shortest cycles containing $e_1$ must also contain at least one of $\VEC e2t$.
		By Lemma \ref{egr: non-incident edges} only $(t-1)(k-1)^{\CFR{g-5}2}$ of them can be accomodated.
	\end{proof}
\end{corollary}

\begin{lemma}\label{egr: counting Bt}
	Let $G$ be an $egr(v,k,g,\lambda)$-graph, with $g=2t$ even.
	For an arbitrary, but fixed vertex $v_0$, let $n_i$ the number of vertices in $B_t=B_t(v_0)$ that has exactly $i$ edges to $B_{t-1}$.
	Note that $n_0=0$ by definition of $B_t$.
	Then following identities hold:
	\begin{enumerate}
		\item $\SE i1kin_i=k(k-1)^{t-1}$.
		\item $\SE i1k\CH i2n_i=\frac{k\lambda}2$.
	\end{enumerate}
	\begin{proof}
		The first identity holds, because both sides count all edges between $B_{t-1}$ and $B_t$.
		$B_{t-1}$ has $k(k-1)^{t-2}$ vertices, and from each of them $k-1$ edges go to $B_t$.
		The lhs counts all edges arriving at $B_t$ from $B_{t-1}$.
		
		For the second identity, we recall that $\frac{k\lambda}2$ shortest cycles contain $v_0$ (lemma \ref{egr: counting cycles using a vertex}). Each such cycle uses exactly one vertex of $B_t$ and is completely contained in $D_t$.
		Conversely every pair of edges between $B_t$ and $B_{t-1}$ with a common endpoint on $B_t$ defines exactly one shortest cycle containing $v_0$. A vertex of $B_t$ that is endpoint of $i$ edges to $B_{t-1}$ is part of $\CH i2$ such pairs, and this proves the second identity.
	\end{proof}
\end{lemma}

\begin{lemma}\label{egr3: shortest cycles on P3}
	Let $G$ be an $egr(v,3,g,2t)$-graph.
	Then every subgraph $P_3$ in $G$ is contained in exactly $t$ shortest cycles.
	\begin{proof}
		Let $P=pxq$ be a $P_3$ in $G$. Let $y$ be the neighbor of $x$ that is not on $P$.
		Each of the $2t$ shortest cycles containing $px$ must contain either $xy$ or $xq$.
		Suppose at least $t+1$ of those cycles contain $xy$.
		Then at most $t-1$ of them contain $xq$.
		This means that there are at least $t+1$ other shortest cycles containing $xq$, and none of them can contain $px$, so they must all contain $xy$. This means that $xy$ is on at least $2t+2$ shortest cycles. Contradiction.
		So at most $t$ of those shortest cycles contain $xy$. By symmetry the same holds for $xq$, which proves the claim.
	\end{proof}
\end{lemma}

\begin{lemma}\label{egr3: connectivity 1}
	Let $G$ be a connected $egr(v,3,g,\lambda)$-graph.
	When $\lambda>2^{\CFR{g-5}2}$ then $G$ is $3$-edge-connected.
	\begin{proof}
		$G$ is clearly $2$-edge-connected, since a bridge cannot be on any cycle.
		Suppose $e,f$ are $2$ edges whose removal disconnects $G$.
		
		If $e$ and $f$ have a common vertex $v$, then the third edge incident with $v$ will disconnect $G$ as well, contradicting that $G$ is $2$-edge-connected.
		So $e$ and $f$ are non-incident edges.
		Since every shortest cycle that contains $e$ must contain $f$ as well, $e$ and $f$ are contained in $\lambda$ different shortest cycles, contradicting Lemma \ref{egr: non-incident edges}.
	\end{proof}
\end{lemma}

\section{$egr(v,4,4,4)$-graphs.}

\begin{figure}[h]
	\centering
	\includegraphics[scale=0.65]{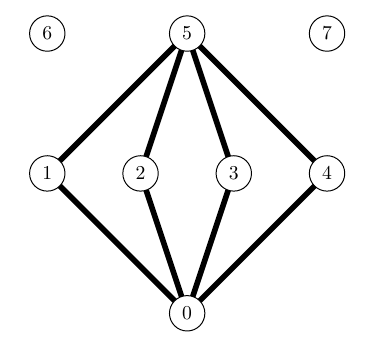}
	\caption{Developing $D_2$ for an $egr(v,4,4,4)$-graph, case 1}\label{egr444c: picture}
\end{figure}

\begin{figure}[h]
	\centering
	\includegraphics[scale=0.65]{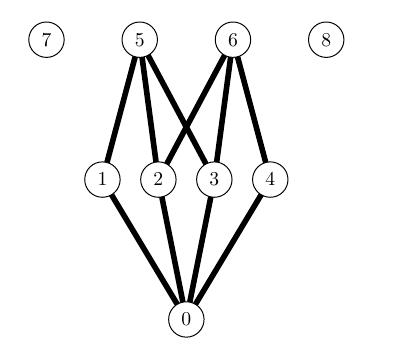}
	\caption{Developing $D_2$ for an $egr(v,4,4,4)$-graph, case 2}\label{egr444b: picture}
\end{figure}

\begin{theorem}
	An $egr(v,4,4,4)$-graph does not exist.
	\begin{proof}
		Suppose $G$ is an $egr(v,4,4,4)$-graph.
		We will take an arbitrary (but fixed) vertex $v_0$ of $G$ and focus on what happens around $v_0$.
		We will write $D_d$ and $B_d$ instead of $D_d(v_0)$ and $B_d(v_0)$.
		We will label $v_0$ with $0$ and use vertex labels in the proof.
		The neighbors of $v_0$ get vertex labels $1$, $2$, $3$, and $4$.
		
		Using the notation of lemma \ref{egr: counting Bt}, the equations become\\
		$4n_4+3n_3+2n_2+n_1=4\cdot3=12$ and\\
		$6n_4+3n_3+n_2=\frac{4\cdot4}2=8$.\\
		
		Case 1: $n_4$ is positive.
		The second equation shows that $n_4=1$, $n_3=0$, and $n_2=2$.
		Substituting in the first equation gives $n_1=4$.
		We will label the vertex of $B_2$ with $4$ edges to $B_1$ with $5$.
		The $2$ vertices of $B_2$ that have $2$ edges to $B_1$ get labels $6$ and $7$, and there are $4$ vertices of $B_2$ with exactly one edge to $B_1$ (part of this is shown in picture \ref{egr444c: picture}).
		By symmetry we may assume that the neighbors of vertex $6$ in $B_1$ are $1$ and $2$.
		Let $x$ be the $4$th neighbor of vertex $1$.
		At this point edge $01$ is contained in $4$ shortest cycles ($0152$, $0153$, $0154$, and $0162$) and edge $15$ is contained in $4$ shortest cycles ($1520$, $1530$, $1540$, and $1526$), and none of them involves vertex $x$.
		So the $4$ shortest cycles containing edge $1x$ all must contain edge $16$ as well, contradicting Lemma \ref{egr3: shortest cycles on P3}.
		
		Case 2: $n_4=0$.
		Subtracting the resulting equations gives $n_1+n_2=4$.
		The first equation reduces to $3n_3+n_2=8$.
		This forces $n_3\leq2$.
		If $n_3=0$, we get $n_2=8$, which violates the first equation.
		If $n_3=1$, we get $n_2=5$, which violates the first equation.
		So $n_3=2$, $n_2=2$, and $n_1=2$.
		We will label the two vertices of $B_2$ with $3$ edges to $B_1$ with $5$ and $6$.
		The two vertices of $B_2$ with $2$ edges to $B_1$ get labels $7$ and $8$, and there are $2$ vertices of $B_2$ with exactly one edge to $B_1$.
		
		Suppose vertices $5$ and $6$ have the same set of neighbors in $B_1$. By symmetry we may assume these neighbors are $1$, $2$, and $3$. Then the vertices $\{1,2,3,0,5,6\}$ induce $K_{3,3}$, which is an $egr(6,3,4,4)$, contradicting Lemma \ref{egr: forbid certain egr subgraphs}.
		
		By symmetry we may assume edges $15,25,35,26,36,46$. At this point the graph looks like Picture \ref{egr444b: picture}.
		Assume vertex $x$ is the $4$th neighbor of vertex $2$.
		At this point edge $02$ is contained in $4$ shortest cycles, and $25$ and $26$ are both contained in $3$ shortest cycles, and none of them involves vertex $x$.
		So at most $2$ of the $4$ shortest cycles containing edge $2x$ can be accommodated. This final contradiction finishes the proof.
	\end{proof}
\end{theorem}

\section{$egr(v,3,8,2t)$-graphs.}

For easy reference we summarize a few of the general lemmas for the case that $G$ is an $egr(v,3,8,2t)$-graph.

\begin{lemma}\label{egr38x: intersection of cycles}
	~
	
	\begin{enumerate}
		\item Every subgraph $P_6$ in $G$ is contained in at most $1$ shortest cycle.
		\item Every subgraph $P_5$ in $G$ is contained in at most $2$ shortest cycles.
		\item Every subgraph $P_4$ in $G$ is contained in at most $4$ shortest cycles.
		\item Every subgraph $P_3$ in $G$ is contained in exactly $t$ shortest cycles.
		\item $G$ is $3$-edge-connected when $t>2$.
		\item Removal of $3$ pairwise non-incident edges cannot disconnect $G$.
	\end{enumerate}
	\begin{proof}
		The claims are special cases of Lemma \ref{egr: intersection of cycles 2}, Lemma \ref{egr3: shortest cycles on P3}, Lemma \ref{egr3: connectivity 1}, and Lemma \ref{egr: non-incident edges}.
	\end{proof}
\end{lemma}

\begin{figure}[h]
	\centering
	\includegraphics[scale=0.65]{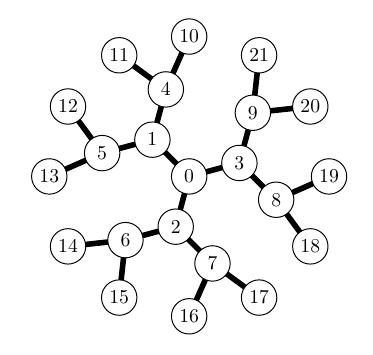}
	\caption{A spanning tree for $D_3$}\label{egr3814: picture}
\end{figure}

Let $G$ be an $egr(v,3,8,2t)$-graph.
We aim to prove non-existence for certain values of $t$, and since every connected component of $G$ must be an $egr(v,3,8,2t)$-graph itself, we may assume that $G$ is connected.

We will take an arbitrary (but fixed) vertex $v_0$ of $G$ and focus on what happens around $v_0$.
We will write $D_d$ and $B_d$ instead of $D_d(v_0)$ and $B_d(v_0)$.
$D_3$ has a spanning tree that looks like Figure \ref{egr3814: picture}, where $v_0$ is the central vertex (labeled with $0$).
Because there are no cycles of length less than $8$, any additional edges in $D_4$ are either in $B_4$ or between $B_3$ and $B_4$.
We see that $D_3$ has exactly $22$ vertices, and $B_3$ has exactly $12$ vertices (labeled $10,\ldots,21$ in the picture).

The $3$ subtrees rooted at the neighbors of $v_0$ ($1,2,3$ in the picture) will be called \emph{branches} and they are labeled by their root (so in the picture branch $1$ contains the vertices $1,4,5,10,11,12,13$).
The intersection of a branch with $B_3$ will be called a branch of $B_3$.

\section{$egr(v,3,8,14)$-graphs.}

\begin{theorem}
	An $egr(v,3,8,14)$-graph does not exist.
	\begin{proof}
		Let $G$ be an $egr(v,3,8,14)$-graph.
		Using the notation of lemma \ref{egr: counting Bt}, we find that the equations become\\
		$3n_3+2n_2+n_1=3\cdot2^3=24$ and\\
		$3n_3+n_2=\frac{3\cdot14}2=21$.\\
		From the second equation we get that $n_2$ must be a multiple of $3$, so by the first equation $n_1$ must also be a multiple of $3$.
		Subtracting the equations we get $n_2+n_1=3$. This reduces the problem to just two cases.
		
		Case 1: $n_1=3$, $n_2=0$, $n_3=7$.
		We refer to picture \ref{egr3814: picture} and focus on branch $1$.
		There are $14$ shortest cycles containing edge $01$. By lemma \ref{egr3: shortest cycles on P3}, $7$ of them must use edge $14$ and $7$ must use edge $15$.
		From Lemma \ref{egr38x: intersection of cycles}.3, and symmetry, it follows that we may assume that from the $7$ shortest cycles that use both $01$ and $14$, $4$ will use edge $4,10$ and $3$ will use edge $4,11$.
		Then from Lemma \ref{egr38x: intersection of cycles}.3 it follows that every edge from $B_3$ to $B_4$ in branch $1$ is contained in at least $1$ of the shortest cycles containing edge $01$.
		
		This implies that every edge from $B_3$ to $B_4$ goes to a vertex of $B_4$ that has at least $2$ edges to $B_3$, so $n_1=0$.
		Contradiction.

		Case 2: $n_1=0$, $n_2=3$, $n_3=6$.
		All vertices in $D_4$ have degree $3$ (in $D_4$) except the $3$ vertices of $B_4$ that have only $2$ edges to $B_3$.
		If there is an edge in $B_4$, we  are left with a single vertex in $D_4$ of degree $2$, which would exhibit a bridge in $G$. Contradiction.

		So there are $3$ edges $e_1$, $e_2$ and $e_3$ "sticking out" of $D_4$, and they form an edge cut of $G$.
		The endpoints of these edges in $D_4$ are all different by construction.

		If these $3$ edges have a common endpoint $u$, every vertex in $D_4$ has degree $3$ in $D_4$, so $D_4$ must be all of $G$. We have $22$ vertices in $D_3$, $6$ vertices of degree $3$ in $D_4$, $3$ vertices of degree $2$ in $D_4$, and $u$, for a total of $32$ vertices. This is a contradiction, since the results in \cite{GJ1} show that there can be no $egr(v,3,8,14)$-graph on $32$ vertices.

		If $2$ of these edges, say $e_1$ and $e_2$, have a common endpoint $u$, let $e_4$ be the third edge incident with $u$.
		Then $\{e_3,e_4\}$ is an edge cut of size $2$, contradicting Lemma \ref{egr38x: intersection of cycles}.5.

		So $e_1$, $e_2$ and $e_3$ are pairwise non-incident.
		Each of the $14$ shortest cycles containing $e_1$ must also contain either $e_2$ or $e_3$, but by Lemma \ref{egr38x: intersection of cycles}.5 only $4+4=8$ can be accommodated. Contradiction.
	\end{proof}
\end{theorem}

\section{$egr(v,k,g,(k-1)^{\frac{g-1}2}-1)$-graphs for $g$ odd and $g\geq5$.}

\begin{figure}[h]
	\centering
	\includegraphics[scale=0.65]{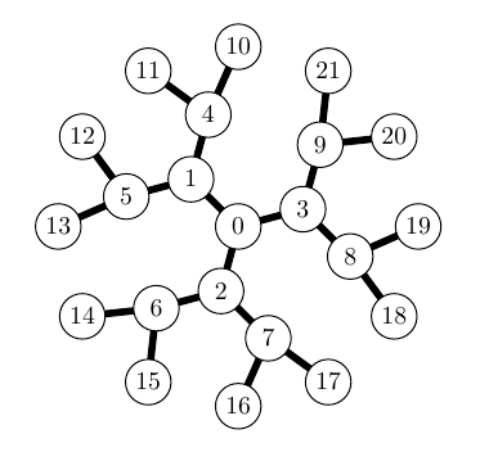}
	\caption{A spanning tree for $D$}\label{egr376: picture}
\end{figure}

\begin{lemma}\label{egr: restrictions near the upper limit}
	Let $G$ be an $egr(v,k,g,(k-1)^{\frac{g-1}2}-1)$-graph, with $g$ odd, and $g\geq5$.
	Then $G$ has exactly $1+k+k(k-1)+\ldots+k(k-1)^\frac{g-3}2+1=2+k\frac{(k-1)^{\frac{g-1}2}-1}{k-2}$ vertices.
	\begin{proof}
		Let $t=\frac{g-1}2$, and $\lambda=(k-1)^t-1$.
		We focus on an arbitrary, but fixed vertex $v_0$ of $G$.
		Let $D=D_t(v_0)$ and $B=B_t(v_0)$ (see \ref{definitions}).
		Picture \ref{egr376: picture} shows a spanning tree for $D$ for the case that $k=3$ and $g=7$.
		We will use a spanning tree $T$ of this type, with $v_0$ being the central vertex, labeled $0$.
		$D$ can (and will) have additional edges, but because of the girth restriction these are all in the outer layer, i.e. in $B$.
		The neighbors of $v_0$ have labels $1,\ldots,k$.
		The subtree with root $i$ will be called \emph{branch} $i$ for $i\in[k]$.
		
		All edges in $B$ must be between vertices of different branches (or we find a cycle of length less than $g$), so every edge in $B$ is on a shortest cycle containing $v_0$ (walk to $v_0$ from each endpoint of the edge).
		
		Conversely, every shortest cycle containing $v_0$ uses exactly one edge of $B$: indeed, the path from $v_0$ to a leaf vertex and then using an edge of $B$ to cross to another branch is a $P_{t+2}=P_{\frac{g+3}2}$, and Lemma \ref{egr: intersection of cycles} tells us that this path is contained in at most one shortest cycle.
		
		By Lemma \ref{egr: counting cycles using a vertex} this means that, just to complete the shortest cycles containing $v_0$, there must be at least $\frac{k\lambda}2=\frac{k(k-1)^t}2-\frac k2$ edges in $B$.
		There cannot be any other edge in $B$, since each additional edge creates an additional shortest cycle containing $v_0$, and we already accounted for all of them.
		
		$B$ has $k(k-1)^{t-1}$ vertices, and every vertex of $B$ is incident with exactly one leaf edge of $T$ that is not contained in $B$, so the maximum degree in $B$ of vertices in $B$ is at most $k-1$. If there were no vertices outside $D$, $B$ would need to have $\frac{k(k-1)^{t-1}\cdot(k-1)}2=\frac{k(k-1)^t}2$ edges. So we are missing $\frac k2$ edges, and this implies there are exactly $k$ edges "sticking out" of $D$. We will call these edges \emph{bad} edges. The endpoint of a bad edge that is not in $D$ will be called a bad vertex.
		
		Branch $1$ must accomodate $\lambda=(k-1)^t-1$ shortest cycles containing edge $\{0,1\}$, and we have seen that each such cycle uses exactly one edge in $B$ with one endpoint in branch $1$ (because all edges in $B$ must be between vertices of different branches). Since branch $1$ contains $(k-1)^{t-1}$ vertices of $B$, this implies that branch $1$ contains exactly one vertex of degree $k-2$ in $B$, and all other vertices that are both in $B$ and in branch $1$ have degree $k-1$ in $B$. By symmetry a similar statement holds for every other branch, so every branch contains exactly one endpoint of a bad edge.
		
		The $k$ bad edges $\VEC e1k$ form an edge cut, so each of the $\lambda$ shortest cycles containing $e_1$ must also use at least one of the edges $\VEC e2k$.
		Applying Lemma \ref{egr: intersection of cycles 2} (with $u-i=3$, so $i=u-3=\frac{g-3}2=t-1$), we see that for an edge $e_i$ that is incident with $e_1$, at most $(k-1)^{t-1}$ edges can contain both $e_1$ and $e_i$.
		If $e_i$ is not incident with $e_1$, we use Lemma \ref{egr: non-incident edges} and find that at most $(k-1)^{t-2}$ shortest cycles containing $e_1$ can also contain $e_i$.
		Even if there is only one edge $e_i$ that is non-incident with $e_1$ we can accomodate at most $x=(k-2)(k-1)^{t-1}+(k-1)^{t-2}$ shortest cycles. We want to show $x<\lambda=(k-1)^t-1$, which finishes the proof, since it shows that all the $e_i$ must be incident with a common vertex of degree $k$, and hence we found all vertices of $G$.
		
		Rearranging the inequality gives $((k-1)^2-(k-2)(k-1)-1)(k-1)^{t-2}>1$, or $(k-2)(k-1)^{t-2}>1$.
		Since $g\geq5$ (so $t\geq2$) this is always true.
	\end{proof}
\end{lemma}

\begin{corollary}
	An $egr(v,4,5,8)$-graph does not exist.
	An $egr(v,6,5,24)$-graph does not exist.
	\begin{proof}
		If an $egr(v,4,5,8)$-graph exists, Lemma \ref{egr: restrictions near the upper limit} tells us that it has exactly $18$ vertices and the results of \cite{GJ1} tell us that this is impossible.
		\\
		Similarly, if an $egr(v,6,5,24)$-graph exists, Lemma \ref{egr: restrictions near the upper limit} tells us that it has exactly $38$ vertices and the results of \cite{GJ1} tell us that this is impossible.
	\end{proof}
\end{corollary}

\section{$egr(v,3,7,6)$-graphs.}

In this section we will assume that $G$ is a $egr(v,3,7,6)$-graph.
Our goal is to prove that such graph does not exist.
Since every connected component of $G$ must be an $egr(v,3,7,6)$-graph itself, we may assume that $G$ is connected.
For easy reference we summarize a few of the general lemmas for this case:

\begin{lemma}\label{egr376: intersection of cycles}
	~
	
	\begin{enumerate}
		\item Every subgraph $P_5$ in $G$ is contained in at most $1$ shortest cycle.
		\item Every subgraph $P_4$ in $G$ is contained in at most $2$ shortest cycles.
		\item Every subgraph $P_3$ in $G$ is contained in exactly $3$ shortest cycles.
		\item $G$ is $3$-edge-connected.
		\item Removal of $3$ pairwise non-incident edges cannot disconnect $G$.
	\end{enumerate}
	\begin{proof}
		The claims are special cases of Lemma \ref{egr: intersection of cycles 2}, Lemma \ref{egr3: shortest cycles on P3}, Lemma \ref{egr3: connectivity 1}, and Corollary \ref{egr: connectivity 2};
	\end{proof}
\end{lemma}

We will take an arbitrary (but fixed) vertex $v_0$ of $G$ and focus on what happens around $v_0$.
We will abbreviate $D_3(v_0)$ as $D$  and $B_3(v_0)$ as $B$.

$D$ has a spanning tree that looks like figure \ref{egr376: picture}, where $v_0$ is the central vertex (labeled with $0$).
$D$ has additional edges but all of those are entirely inside $B$ (otherwise we find cycles with length less than $7$).
We see that $D$ has exactly $22$ vertices, and $B$ has exactly $12$ vertices (labeled $10,\ldots,21$ in the picture).

The three subtrees rooted at the neighbors of $0$ ($1,2,3$ in the picture) will be called branches and they are labeled by their root (so in the picture branch $1$ contains the vertices $1,4,5,10,11,12,13$).
The intersection of a branch with $B$ will be called a branch of $B$.
The $2$ pairs in a branch of $B$ that share a common neighbor will be call subbranches (in the picture $\{10,11\}$ and $\{12,13\})$ are subbrances of branch $1$ of $B$.

\begin{lemma}\label{egr376: even and odd vertices}
	Let $u$ and $w$ vertices from the same subbranch of $B$.
	Then one of them has degree $1$ in $B$, and one of them has degree $2$ in $B$.
	\begin{proof}
		Let $x$ be the unique common neighbor of $u$ and $w$, and $y$ the unique common neighbour of $v_0$ and $y$ (in the picture: if we take $w=10$, then $u$ must be $11$, $x$ must be $4$, and $y$ must be $1$).
		Then $P=wxyv_0$ is the unique shortest $w,v_0$-path in $G$.
		By Lemma \ref{egr376: intersection of cycles}.3 there are exactly $3$ shortest cycles containing $v_0yx$ and by Lemma \ref{egr376: intersection of cycles}.2 at most $2$ of them also contain $u$, and at most $2$ of them also contain $w$.
		By symmetry we may assume that $2$ of them contain $w$, so only one contains $u$.
		Then $w$ has degree $2$ in $B$, and $u$ has degree at least $1$ in $B$.
		If $u$ also has degree $2$ in $B$, then there are $2$ shortest cycles containing $v_0yxu$: first walk $uxyv_0$, then a shortest path to any neighbor of $u$, and then back to $u$. Since $u$ and $w$ have no common neighbor, we find $4$ different shortest cycles containing the $P_3$ $v_0,y,x$, contradicting Lemma \ref{egr376: intersection of cycles}.3.
	\end{proof}
\end{lemma}

So in each subbranch there is one vertex of degree $1$ in $B$, which we will call an odd vertex, and one vertex of degree $2$ in $B$, which we will call an even vertex.
Referring to the picture this means that, by appropriate relabeling, we may assume that the even vertices correspond to even labels, and the odd vertices correspond to odd labels.

\begin{lemma}\label{egr376: paths in B}
	Following statements are true:
	\begin{enumerate}
		\item $B$ cannot contain a cycle.
		\item $B$ cannot contain a path of length longer than $7$.
		\item $B$ contains exactly $3$ paths, all of positive length.
		\item An edge in $B$ cannot connect $2$ vertices of the same branch.
		\item A path in $B$ cyclically visits the cycles (in the picture: if a $P_5$ start in branch $1$ and goes counterclockwise, it will consecutively visit the branches $12312$).
		\item Consecutive revisits of the same branch hit different subbranches (with the example of the previous item: if the $P_5$ starts at vertex $11$, the next visit to branch $1$ will be at vertex $12$; note that the first vertex of a path in $B$ is always odd, and for a $P_5$ the fourth vertex is even).
	\end{enumerate}
	\begin{proof}
		~
		
		\begin{enumerate}
			\item $B$ cannot contain a cycle.\\
			The shortest possible cycle already has $7$ vertices, and Lemma \ref{egr376: even and odd vertices} tells use that among the $12$ vertices of $B$ there are only $6$ vertices of degree $2$.
			
			\item $B$ cannot contain a path of length longer than $7$.\\
			Such path contains at least $7$ vertices of degree $2$.
			
			\item $B$ contains exactly $3$ paths, all of positive length.\\
			$B$ has no cycle, $12$ vertices, $6$ vertices of degree $1$, and $6$ vertices of degree $2$.
			
			\item An edge in $B$ cannot connect $2$ vertices of the same branch.\\
			Such an edge would create a cycle of length at most $5$.
			
			\item A path in $B$ cyclically visits the cycles.\\
			If a path would revisit the same branch after $2$ steps, this creates a cycle of length at most $6$.
			 (In the picture: if $B$ contains path $10,14,12$, we find $C_6$ $10,14,12,5,1,4$.)
			
			\item Consecutive revisits of the same branch hit different subbranches.\\
			A path revisits the same branch after three steps. If both visits would hit the same subbranch we find a cycle of length $5$. (In the picture: if $B$ contains path $10,14,18,11$, we find $C_5$ $10,14,18,11,4$.)
		\end{enumerate}
	\end{proof}
\end{lemma}

Since every odd vertex of $B$ has degree $1$ in $B$ (and degree $2$ in $D$), each of them is incident with exactly one edge whose other endpoint is not in $D$. We will call those edges \emph{bad edges} and label them $f_{11},f_{13},\ldots,f_{21}$ (to indicate to which vertex they are attached).
The endpoints of bad edges that are not in $D$ will be called \emph{bad vertices}.

Note that the bad edges form an edge cut in $G$.
Also note that the bad edges attached to vertices of the same branch of $B$ cannot have a common endpoint, since that would produce a cycle of length at most $6$. (In the picture: if $f_{11}$ and $f_{13}$ have a common endpoint $x$, we find a $C_6$ with vertices $x,13,5,1,4,11$.)

\begin{lemma}\label{egr376: bad edges}
	Following claims are true:
	\begin{enumerate}
		\item $3$ bad edges cannot have a common endpoint.
		\item There cannot be $3$ sets of $2$ bad edges with a common endpoint.
	\end{enumerate}
	\begin{proof}
		~
		
		\begin{enumerate}
			\item $3$ bad edges cannot have a common endpoint.\\
			If $x$ is the common endpoint of $3$ bad edges, the $3$ remaining bad edges still form an edge cut. If $2$ of them are incident at vertex $w$, we replace them by the third edge incident with $w$, and we found an edge cut of size $2$, contradicting Lemma \ref{egr376: intersection of cycles}.4.
			If they are pairwise non-incident we contradict Lemma \ref{egr376: intersection of cycles}.5.
			
			\item There cannot be $3$ sets of $2$ bad edges with a common endpoint.\\
			We replace each pair by the third edge incident with the common endpoint, and again find an edge cut of size $3$.
			If these $3$ edges have a common endpoint, we have found all ($22+3+1=26$) vertices and all edges of $G$ and the results of \cite{GJ1} prove that no $egr(26,3,7,6)$ exists. At this point we can finish the proof as in item 1.
		\end{enumerate}
	\end{proof}
\end{lemma}

\begin{theorem}\label{egr376: main theorem}
	Let $w$ be an even vertex of $B$, $u_e$ the even vertex in the other subbranch of the branch containing $w$, and $u_o$ the odd vertex in the same subbranch as $u_e$.
	Then exactly one of the following conditions is true:
	\begin{enumerate}
		\item $w$ and $u_e$ are on a path in $B$ of length at least $5$ and they have distance $3$ on that path.
		\item $w$ and $u_o$ are on a path in $B$ of length at least $4$ and they have distance $3$ on that path.
		\item The bad edge at $u_o$ has a common neighbor with the bad edge at a neighbor of $w$ (and this implies that $w$ must be a penultimate vertex on the maximal path containing $w$ in $B$).
	\end{enumerate}
	\begin{proof}
		Let $P=wxyv_0$ be the unique shortest $w,v_0$-path.
		Let $z$ be the neighbor of $y$ that has no name yet, so $u_e$ is the even neighbor of $z$ in $B$ and $u_o$ is the odd neighbor of $z$ in $B$.
		(Picture reference: if $w=10$, then $x=4$, $y=1$ $z=5$, $u_e=12$ and $u_o=13$.)
		\\
		By Lemma \ref{egr376: intersection of cycles}.3 exactly $3$ shortest cycles contain the path $wxy$.
		$2$ of them contain $v_0$: the path $v_0yxw$ can be extended with each neighbor of $w$ ($2$ such neighbors exist, since $w$ is an even vertex) and then completed to a shortest cycle.
		By Lemma \ref{egr376: intersection of cycles}.2 the remaining shortest cycle $C$ cannot contain $v_0$, and hence must contain $z$, and either $u_e$ or $u_o$.
		So far we found that $C$ contains the path $Q=wxyzu$ where $u$ is either $u_e$ or $u_o$.
		At the $w$-side we can only prolong $Q$ with a neighbor $w'$ of $w$ in all cases.
		
		If the $u=u_e$ we only can prolong $Q$ with a neighbor $u'$ of $u_e$.
		Now we have a path of $7$ vertices $w'wxuzu_eu'$ which must be completed to a $C_7$, so $w'u'$ is an edge of $B$.
		This is case 1: $ww'u'u_e$ is a path $R$ of length $3$, and because both $u_e$ and $w$ are even vertices, $R$ can be prolonged on both sides to a path of length at least $5$.
		
		If $u=u_o$ we have two possibilities. Either we prolong $Q$ (at the $u_o$-side) with a vertex $u'$ along an edge in $B$.
		Again we have a path of $7$ vertices $w'wxuzu_eu'$ which must be completed to a $C_7$, so $w'u'$ is an edge of $B$.
		This is case 2: $ww'u'u_e$ is a path $R$ of length $3$, and because $w$ is an even vertex, $R$ can be prolonged on one side to a path of length at least $4$.
		
		The last possibility is, that we continue at the $u_o$-side using the bad edge at $u_o$, say to vertex $p$.
		In order to complete a $C_7$, $pw'$ must be an edge of $G$.
		This is only possible if $w'$ is an odd vertex (even vertices have no neighbors outside $D$), so $w$ must have been a penultimate vertex, and $pw'$ is the bad edge at $w'$. This is exactly case $3$.
		
		Finally, each condition defines a different $C_7$: case 1 uses $u_e$ but not $u_o$, case $2$ uses $u_o$ but not $u_e$ and case $3$ is the only case that uses a vertex outside $D$. If more than one condition is true we find more than one $C_7$ containing $wxyz$ and hence more than $3$ $7$-cycles using $wxy$, contradicting Lemma \ref{egr376: intersection of cycles}.3.
	\end{proof}
\end{theorem}

\begin{corollary}\label{egr376: no P5}
	$B$ cannot contain $P_5$ as a maximal path.
	\begin{proof}
		Let $w$ be the central vertex of a maximal path $P$ of length $4$. Since $w$ is not penultimate, case 3 cannot occur.
		Since there are no vertices on $P$ at distance $3$ from $w$ the first two cases cannot occur either.
	\end{proof}
\end{corollary}

For the remainder of the argument we will choose an endpoint of a path as its first vertex in such a way that walking the path from the first vertex to the last will hit the branches counterclockwise (referring to Picture \ref{egr376: picture}), i.e. in the cyclic repetition $1,2,3,1,2,3,\ldots$.
Also, for a path $P$ we will label the $i$-th vertex as $P(i)$ (starting at index $1$).
We will also label the subbranches of branch $1$ with $A_1,A_2$, branch $2$ with $B_1,B_2$, and branch $3$ with $C_1,C_2$.
Because these subbranch names are all different we need not specify the branch number when we know the subbranch name.

\begin{lemma}\label{egr376: no P8}
	$B$ cannot contain $P_8$.
	\begin{proof}
		Any $P_8$ in $B$ must be a maximal path by Lemma \ref{egr376: paths in B}.2.
		Suppose $P$ is a $P_8$ in $B$.
		By symmetry we may assume that $P(1)$ is the odd vertex in subbranch $A_1$, and $P(2)$ is the even vertex in subbranch $B_1$.
		By Lemma \ref{egr376: paths in B}.6 revisiting a branch must hit a different subbranch, so $P(4)$ must be the even vertex in subbranch $A_2$, and $P(5)$ must be the even vertex in subbranch $B_2$.
		Then $P(7)$ (an even vertex) is in subbranch $A_1$ again and $P(8)$ (an odd vertex) is in subbranch $B_1$ again.
		This produces $2$ edges between subbranch $A_1$ and subbranch $B_1$, the first one is edge $P(1)P(2)$, the second one is edge $P(7)P(8)$. Combining these with the common neighbors of those subbranches we find $C_6$. Contradiction.
		(Picture reference: we may assume $P$ consecutively visits labels $11,14,18,12,16,20,10,15$ (after choosing the first $3$ vertices of $P$, the rest of the path is forced), and we find $6$-cycle $4,10,15,6,14,11$).
	\end{proof}
\end{lemma}

\begin{lemma}\label{egr376: no P7}
	$B$ cannot contain $P_7$.
	\begin{proof}
		Lemma \ref{egr376: no P8} excludes $P_8$, so any $P_7$ in $B$ must be a maximal path.
		Suppose $P$ is a $P_7$ in $B$.
		We may assume that $P(1)$ is the odd vertex in subbranch $A_1$.
		By Lemma \ref{egr376: paths in B}.6 revisiting a branch must hit a different subbranch.
		So $P(4)$ visits subbranch $A_2$, and $P(7)$ must be the odd vertex in subbranch $A_1$ again, contradicting that $P$ is a path.
		(Picture reference: we may assume $P$ consecutively visits labels $11,14,18,12,16,20,11$.)
	\end{proof}
\end{lemma}

\begin{lemma}\label{egr376: no P6}
	$B$ cannot contain $P_6$.
	\begin{proof}
		Suppose $P$ is a $P_6$ in $B$. As in the previous lemmas, $P$ must be a maximal path.
		By symmetry, we may assume that $P(1)$ is the odd vertex in $A_1$, and $P(3)$ is the even vertex in $C_1$.
		By Lemma \ref{egr376: paths in B}.6 revisiting a branch must hit a different subbranch.
		So $P(4)$ must be in $A_2$, which means that $P(3)P(4)$ connects the even vertex in $C_1$ to the even vertex in $A_2$ (*).
		$P$ occupies an odd vertex in branch $1$ and branch $3$. Since $P(3)$ is in $C_1$, $P(6)$ must be the odd vertex in $C_2$.
		
		$P$ has $4$ even vertices, so there are only two options for the other $2$ paths in $B$:
		\begin{enumerate}
			\item One is $P_4$ and the other is $P_2$.
			\item They are both $P_3$.
		\end{enumerate}
		
		We consider the first case.
		Let $Q_1$ be the $P_4$ in $B$ and $Q_2$ the $P_2$ in $B$. $Q_1$ must start and end in the same branch, i.e. it exhausts the odd vertices of a branch. Since $P$ occupies an odd vertex in both branch $1$ and branch $3$, $Q_1$ must start and end in branch $2$.
		This forces $Q_2$ to join the remaining odd vertex in branch $3$ to the remaining odd vertex in branch $1$.
		Since $P(6)$ occupies the odd vertex in $C_2$, $Q_2$ must use the odd vertex in $C_1$.
		Since $P(1)$ occupies the odd vertex in $A_1$, $Q_2$ must use the odd vertex in $A_2$.
		So $Q_2$ connects the odd vertex in $C_1$ to the odd vertex in $A_2$ (**).
		\\
		(*) and (**) combine with the common neighbors of $C_1$ and $A_2$ to a $C_6$. Contradiction.
		(Picture reference: assuming $A_1=\{10,11\})$ and $C_1=\{18,19\}$, we get $P(4)=12,P(3)=18,Q_2(1)=19,Q_2(2)=13$ and we find the $6$-cycle $13,5,12,18,8,19$.)
		
		We consider the second case. Let $Q_2$ and $Q_3$ be the $2$ $P_3$s.
		$P$ has exhausted the even vertices in branch $2$ and since $Q_1$ and $Q_2$ both need one even vertex we may assume that $Q_2$ starts in branch $2$ (so its even vertex $Q_2(2)$ is in branch $3$), and $Q_3$ starts in branch $3$ (so its even vertex $Q_3(2)$ is in branch $1$).
		
		We check branch $1$. It contains the even vertices $P(4)$ and $Q_3(2)$, and the odd vertices $P(1)$ and $Q_2(3)$.
		Since $P(4)$ and $P(1)$ cannot be in the same subbranch (lemma \ref{egr376: paths in B}.6), $P(4)$ (even) shares a subbranch with $Q_2(3)$ (odd), and $Q_3(2)$ (even) with $P(1)$ (odd).
		Applying Theorem \ref{egr376: main theorem} with $w=Q_3(2)$ we find $u_e=P(4)$ and $u_o=Q_2(3)$.
		$w$ and $u_e$ are on different paths, so case 3 of Theorem \ref{egr376: main theorem} applies: we find a common endpoint for the bad edges at $u_o=Q_2(3)$ and a neighbor of $w$, which must be $Q_3(1)$ or $Q_3(3)$.
		Both possibilities give us one pair of incident bad edges, not involving a vertex of $P$.
		
		We check branch $2$. It contains the even vertices $P(2)$ and $P(5)$, and the odd vertices $Q_2(1)$ and $Q_3(3)$, but we do not know how they are combined.
		Applying the theorem with $w=P(2)$ we get a common neighbor for the bad edge at $P(1)$ and either $Q_2(1)$ or $Q_3(3)$.
		Applying the theorem with $w=P(5)$ we get a common neighbor for the bad edge at $P(6)$ and either $Q_2(1)$ or $Q_3(3)$.
		Possibility 1: $P(1)$ has a common neighbor with $Q_2(1)$ and $P(6)$ with $Q_3(3)$.
		Possibility 2: $P(1)$ has a common neighbor with $Q_3(3)$ and $P(6)$ with $Q_2(1)$.
		(The other two possibilities can be excluded because of Lemma \ref{egr376: bad edges}.1.)
		Either possibility gives us $2$ pairs of incident bad edges, both involving a vertex of $P$.
		
		Combining the results of branch 1 and branch 2, we either have $3$ pairs of bad edges, or at least $3$ bad edges with a common endpoint, contradicting Lemma \ref{egr376: bad edges}.
	\end{proof}
\end{lemma}

\begin{theorem}
	There is no such thing as an $egr(v,3,7,6)$-graph.
	\begin{proof}
		In view of all exclusions above, the longest path in $B$ is a $P_4$, which means that $B$ must contain exactly $3$ $P_4$s.
		We'll name them $Q_1,Q_2,Q_3$, where $Q_i$ starts (and ends) in branch $i$.
		In branch $1$ we find as odd vertices the start and end vertex of $Q_1$, and as even vertices the second vertex of $Q_3$ and the third vertex of $Q_2$, but we do not know how they combine in the subbranches, and similar for the other branches.
		Since there are two possible combinations in each branch, this gives us $8$ possibilities.
		Fortunately we can reduce the number of cases.
		If a path start in a given subbranch, the other vertex in that subbranch is occupied by either the third vertex of the $P_4$ that start in the next branch, or by the second vertex of the $P_4$ that starts in the previous branch.
		Using this terminology, we can reduce to two possibilities:
		\begin{enumerate}
			\item In all branches the other vertex is occupied by the third vertex of the cyclically next $P_4$.
			\item In exactly two branches the other vertex is occupied by the third vertex of the cyclically next $P_4$.
		\end{enumerate}
		The other two options can be eliminated by reversing the direction.
		
		We handle the first possibility:
		\\
		In branch 1, $Q_1(1)$ and $Q_2(3)$ are in the same subbranch, so $Q_1(4)$ and $Q_3(2)$ are in the other subbranch.
		Theorem \ref{egr376: main theorem} with $w=Q_2(3)$ gives us a common neighbor for the bad edges at $Q_1(4)$ and $Q_2(4)$.
		\\
		In branch 2, $Q_2(1)$ and $Q_3(3)$ are in the same subbranch, so $Q_2(4)$ and $Q_1(2)$ are in the other subbranch.
		Theorem \ref{egr376: main theorem} with $w=Q_3(3)$ gives us a common neighbor for the bad edges at $Q_2(4)$ and $Q_3(4)$.
		\\
		So the three bad edges at $Q_1(4)$, $Q_2(4)$ and $Q_3(4)$ all have the same bad vertex as endpoint, contradicting Lemma \ref{egr376: bad edges}.1.
		
		Finally, note that this proof does not care at all about what happens in branch $3$, so this proof also covers the second possibility.
	\end{proof}
\end{theorem}

\section{Discussion}

Although all proofs in here are more or less ad hoc, they still share a ''flavor".
It is very probable that similar techniques can be used to increase lower bounds for other $egr(v,k,g,\lambda)$-graphs.
Finding an overarching approach would be great!

\section*{Acknowledgements}

Thanks to Jorik Jooken for carefully reading an early version of this paper, and suggesting many improvements.

\ms
\end{document}